\documentclass[15pt]{article}
\usepackage[paperheight=29.7cm, paperwidth=21cm]{geometry}
\usepackage[T2A,T1]{fontenc}
\usepackage[utf8]{inputenc}
\usepackage{amsmath}
\usepackage{amssymb}
\usepackage{graphicx}
\title{\textbf{BOUNDEDNESS OF THE DİSCRETE HİLBERT TRANSFORM ON DİSCRETE WEIGHTED MORREY SPACES }}
\author{Rashid Aliev$^{1}$, Amil Jabiyev$^{2}$}
\begin{document}
\maketitle
\begin{center}
$^{1,2}$ {Baku State University, Baku, Azerbaijan}\\
\vspace{1mm}
e-mail: $^{1}$aliyevrashid@mail.ru, $^{2}$ amilcbiyev23@gmail.com \end{center}
\textbf{Abstract. } The Hilbert transform is a multiplier operator and is widely used in the theory of Fourier transforms. The Hilbert transform was the motivation for the development of modern harmonic analysis. Its discrete version is also widely used in many areas of science and technology and plays an important role in digital signal processing. The essential motivation behind thinking about discrete transforms is that experimental data are most often not taken in a continuous manner but sampled at discrete time values. Since much of the data collected in both the physical sciences and engineering are discrete, the discrete Hilbert transform is a rather useful tool in these areas for the general analysis of this type of data. In this paper, we discuss the discrete Hilbert transform on discrete Weighted Morrey spaces and obtain its boundedness in these spaces.
\\
\textbf{Keywords:} Discrete Hilbert transform, Morrey spaces, weighted Morrey spaces, discrete Morrey spaces, discrete weighted Morrey spaces, the boundedness of the discrete Hilbert transform on weighted Morrey spaces.

\section{Introduction}
We denote by $l_{p}$, $p \geq1 $, the class of scalar sequences $b= \{b_n\}_{n_\in \mathbb{Z}}$ satisfying the following condition:
\[ {{\|b\|}_{l_p}}=\left({{\sum_{n\in\mathbb{Z}}}{|b_n|}^p}\right)^\frac{1}{p} <\infty ,\]  \\ where $\mathbb{Z}$  is the set of integers. \\ 
\indent Let $b={\{b_n\}}_{n\in\mathbb{Z}} \in l_p , p\geq1$. The sequence $H(b)=\{\ (Hb)_n\} _{n\in \mathbb{Z}}$, where \\
\[\ (Hb)_{n}=\sum_{m\neq n} \frac{b_m}{n-m}, \indent n\in \mathbb{Z},  \] is called the Hilbert transform of the sequence $b= \{b_n\}_{n_\in \mathbb{Z}}$. \\ \indent M.Riesz proved (see \cite{Riesz}) that if $b\in l_p $, $p>1$, then $H(b) \in l_p$ and the inequality
\begin{equation}\label{1.1}
 \| H(b)\|_{l_p} \leq C_p \|b\|_{l_p} 
 \end{equation}
 holds, where $C_p$ is an absolute constant. Weighted analogues of (\ref{1.1}) are investigated in the works \cite{Anders,Belov1,Belov2,Gabis,Hunt,Lifl,Rakot,Step}.
If $b\in l_{1} $ then the sequence $H(b)$ belongs to the class $\bigcap _{p>1}l_{p}  $, but generally does not belong to the class $l_{1} $ (see \cite{Aliev4}). In this case, R.Hunt, B.Muckenhoupt and R.Wheeden (see \cite{Hunt}) proved that the distribution function  $ (Hb)(\lambda ):=\sum _{\{n\in \mathbb{Z}:\, \, |(Hb)_{n} |>\lambda \}}\, \, \, 1 $ of $H(b)$ satisfies the weak condition
\begin{equation*}
\forall \lambda >0  \qquad  |(Hb)(\lambda )|\le \frac{c_{0} }{\lambda } \| b\| _{l_{1} } ,
\end{equation*}
where $ c_{0} $ is an absolute constant. In \cite{Aliev4}, it was proved that, if the sequence $b\in l_{1} $ satisfies the conditions $\sum _{n\in \mathbb{Z}}\, \, \, b_{n}  =0$ (this condition is necessary for the summability of the discrete Hilbert transform) and $\sum _{n\in \mathbb{Z}}\, \, \, |b_{n} |\ln (e+|n|) <\infty $, then $ H(b) \in l_{1} $ and the following inequality holds:
\begin{equation*}
\| H(b) \| _{l_{1} } \le 6\sum _{n\in \mathbb{Z}}\, \, \, |b_{n} |\ln (e+|n|) .
\end{equation*}

In \cite{Aliev3}, the concept of $Q$-summability of series was introduced and using this notion it was proved that the Hilbert transform of a sequence $b\in l_{1} $ is $Q$-summable and its $Q$-sum is equal to zero. In \cite{Aliev1,Aliev2,Aliev3,Aliev5,Hao1,Hao2,Jakfar} discrete analogues of harmonic analysis operators on discrete Morrey and discrete Orlicz spaces were studied.

 In this paper, we discuss the discrete Hilbert transform on discrete Weighted Morrey spaces. In particular, we obtain its boundedness on the discrete Weighted Morrey spaces using the boundedness of the Hilbert transform on Weighted Morrey spaces.
\\
\section{Discrete Weigted Morrey spaces}
The classical Morrey spaces $M_{\lambda ,p}(R) $, $0\le \lambda \le \frac{1}{p} $, $1\le p<\infty $ (see \cite{Adams,Morrey,Sickel1,Sickel2,Sickel3}), consist of the functions $f\in L_{p,loc} \left(R \right)$  for which the following norm is finite
\begin{equation*}
\| f\| _{M _{\lambda ,p} } =\mathop{\sup }\limits_{x} \mathop{\sup }\limits_{r>0} \left[r^{-\lambda } \left\| f\right\| _{L_{p} (B(x;r))} \right].
\end{equation*}

We note that if $\lambda =0$, then $M_{\lambda ,p} =L_{p} $; if $\lambda =\frac{1}{p} $, then $M_{\lambda ,p} =L{}_{\infty } $  (see \cite{Adams}). In case $p>1$, $0\le \lambda <\frac{1}{p} $, F.Chiarenza and M.Frasca (see \cite{Chia}) showed the boundedness of the Hardy--Littlewood maximal operator, the fractional integral operator and a singular integral operator in the Morrey spaces. Hence, in particular, it implies the boundedness of the Hilbert transform in Morrey spaces. It means that, in case $p>1$, $0\le \lambda <\frac{1}{p} $, for any $f(t)\in M_{\lambda ,p} $ we have
\begin{equation*}
(Hf)(t):=\frac{1}{\pi } v.p. \int _{R}\frac{f(\tau)}{t-\tau } d\tau  \in M_{\lambda ,p}, 
\end{equation*}
and there exist $C_{\lambda ,p} >0$ such that
\begin{equation*}
\| Hf\| _{M_{\lambda ,p} } \le C_{\lambda ,p} \cdot \| f\| _{M_{\lambda ,p} } 
\end{equation*}
holds for all $f\in M_{\lambda ,p} $.

In \cite{Gun1}, the authors introduced a discrete analogue of Morrey spaces and studied their inclusion properties. For $m\in Z$ and $n\in N\cup \{0\}$ define $S_{m,n} =\{k\in Z: |k-m|\le n\}$. Following standard conventions, we denote the cardinality of a set $S$ by $|S|$. Then we have $|S_{m,n}|=2n+1$ for all $m\in Z$ and each $n\in N\cup \{0\}$. Discrete Morrey spaces $m_{\lambda ,p} $, $0\le \lambda \le \frac{1}{p} $, $1\le p<\infty $, consist of the sequences $b=\{b_{n} \}_{n\in Z} $  for which the following norm is finite
\begin{equation*}
\| b\| _{m_{\lambda ,p} } =\mathop{\sup }\limits_{m\in Z} \mathop{\sup }\limits_{n\in N\cup \{0\}} \left[|S_{m,n}|^{-\lambda } \left(\sum _{k\in S_{m,n} }|b_{k} |^{p}  \right)^{{1 \mathord{/{\vphantom{1 p}}.\kern-\nulldelimiterspace} p} } \right]. 
\end{equation*}

 In \cite{Gun2}, the boundedness of the discrete Hardy-Littlewood maximal operators and discrete Riesz potentials is proved, and in \cite{Aliev2}, the boundedness of the Hilbert transform on discrete Morrey spaces is proved.
 
Let $0<\lambda<\frac{1}{p}$ , $1\leq p<\infty$ and a weight function $w(t)$ be a locally integrable function on $R$ that takes values in $(0,\infty)$ almost everywhere. The weighted Morrey spaces $M_{\lambda,p,w}$ (see \cite{Komori}), consist functions $f\in L_{p,w,loc}(R)$ for which the following norm is finite 
\\ \[ {\|f\|_{M_{\lambda,p,w}}} = \sup_{x} \sup_{r>0} \frac{\left(\int_{B(x,r)} {|f(t)|^p w(t)dt}\right)^\frac{1}{p} }{\left(\int_{B(x,r)}w(t)dt\right)^\lambda}. \] 
\\
If $w(t) \equiv 1$ , then $M_{\lambda,p,w} = M_{\lambda,p}$, $0<\lambda<\frac{1}{p}$ , $1\leq p<\infty$ . Y.Komori and S.Shirai showed (see \cite{Komori}) the boundedness of the Calderon-Zygmund operator on the weighted Morrey spaces. In particular, this implies the boundedness of the Hilbert transform on the weighted Morrey spaces. This means that, in case $p>1$, $0<\lambda<\frac{1}{p}$, for any $f(t)\in M_{\lambda,p,w}$ , we have 
\\
\[H(f) \in M_{\lambda,p,w} \] and there exists $C_{\lambda,p,w}>0$ such that
\\
\[ ||Hf||_{M_{\lambda,p,w}} \leq C_{\lambda,p,w} ||f||_{M_{\lambda,p,w}}\] holds for all $f\in M_{\lambda,p,w}$ . 
\\

For any $\lambda >0$ and for any interval $B \subset R$ denote by $\lambda B $ the interval with the same center as $B$ whose radius length is $\lambda$ times that of $B$. 
\\  \\
\textbf{Definition 2.1.} (see \cite{Komori}) If there exists a constant $D>0$ such that for any interval $B \subset R$ we have
\[ w(2B) \leq D w(B), \]
then we say that $w$ satisfies the doubling condition and we denote $w\in {\Delta_2}$, where $w(E)=\int_{E} {w(t)dt}$.
\\  \\
\textbf{Definition 2.2.} (see \cite{Hunt})  A weight function $w$ is in the Muckenhoupt class $A_{p}$ with $1<p< \infty$ if there exists $C>0$ such that for any interval $B \subset R$
\begin{equation}\label{2.1}
 \left( \frac{1}{|B|} \int_{B} w(x)dx \right) \left( \frac{1}{|B|} \int_{B} w(x)^{-\frac{1}{p-1}}dx \right)^{p-1} \leq C, 
 \end{equation}
and the infimum of $C$ satisfying the inequality (\ref{2.1}) is denoted by $[w]_{A_{p}}$. We define $A_{\infty}=\cup_{1<p<\infty} A_{p}$.

We will need the following lemmas:
\\  \\
\textbf{Lemma 2.1 (see \cite{Grafakos}).} If $w\in A_p$ for some $1\leq p<\infty$ , then $w\in \Delta_2$.
\\  \\
\textbf{Lemma 2.2 (see \cite{Komori}).} If $w \in \Delta_2$ , then there exists a constant $D_1 >1$ such that for any interval $B$
\[ w(2B)\geq D_1 w(B). \]

For $m=(m_1 ,m_2 ,..., m_d) \in \mathbb{Z}^d$, $n\in \mathbb{N \cup} \{0\}$, we define
\[S_{m,n}:= \left\{ k\in \mathbb{Z}^d :  ||k-m||_\infty \leq n \right\}, \]
where $\|k-m\|_\infty =\max_{1\leq i \leq d} \{ |k_{i}-m_{i} | \} $.

Let $0<\lambda<\frac{1}{p}$ , $1\leq p<\infty$. Discrete Weighted Morrey spaces $m_{\lambda,p,w}$ (see \cite{Hao1}) consists of sequences $b: \mathbb{Z} \to R$ for which the following norm is finite
\[||b||_{m_{\lambda,p,w}} := {\sup_{m\in \mathbb{Z}} \sup_{n\in \mathbb{N}\cup \{0\}}} \frac{\left(\sum_{k\in S_{m,n}}|b_k|^p w_k\right)^\frac{1}{p}}{\left(\sum_{k\in S_{m,n}} w_k\right)^\lambda} <\infty, \]
where $w=\{w_{k}\}_{k \in \mathbb{Z}}$  is weight and $w_k \in (0,\infty)$ for any $k \in \mathbb{Z}^{d}$. If $\{w_k\} \equiv 1$, then $m_{\lambda,p,w} =m_{\lambda,p}$ (see \cite{Gun1}). 
\\  \\
\textbf{Definition 2.3 (see \cite{Hunt}).} A weight $\{w_k\}$ belongs to $\tilde{A}_p$ for $1<p<\infty$ , if it satisfies the following condition:
\[ [w]_{\tilde{A}_p} := \sup_{{m,n}\in \mathbb{Z}} \frac{\left(\sum_{k=m}^{n} w_k\right) \left(\sum_{k=m}^{n} {w_k}^{-\frac{1}{p-1}} \right)^{p-1}}{(n-m+1)^p} <\infty. \]
\\  \\
\textbf{Definition 2.4.} If there exists $D>0$, such that for any $m \in \mathbb{Z}$ and for any $\delta>0$
\begin{equation}\label{2.2}
\sum_{ \{k \in \mathbb{Z}: \ k \in [m-2\delta, m+2\delta] \} } w_k \leq D \sum_{ \{k \in \mathbb{Z}: \ k \in [m-\delta, m+\delta] \} } w_k,
\end{equation}
then we say $\{w_k\} \in \tilde{\Delta}_2$.
\\  \\
\textbf{Lemma 2.3.} If $\{w_k\}\in \tilde{\Delta}_{2}$, then there exists $D_1 >1$ such that for any $m \in \mathbb{Z}$ and $n \in \mathbb{N} \cup \{0\}$
\[ \sum_{k=m-(2n+1)}^{m+(2n+1)} w_k \ge D_1 \sum_{k=m-n}^{m+n} w_k. \] 
\\  \\
\textbf{Lemma 2.4.} If $\{w_k\} \in \tilde{A}_p$ for some $1 \leq p < \infty$, then $ \{w_k\} \in \tilde{\Delta}_2$.
\\  \\
\textit{Proofs of Lemmas 2.3 and 2.4.} To prove the lemmas, we define the function $w(x)=w_k$ for $x \in [k-\frac{1}{2}, k+\frac{1}{2})$, $k \in \mathbb{Z}$. Then for any $k \in \mathbb{Z}$ we have $w_k =\int_{k-\frac{1}{2}}^{k+\frac{1}{2}} w(x)dx$. Let $\{w_k\}\in \tilde{\Delta}_{2}$. Let us prove $w \in \Delta_2$. For this we take any interval $B=[b-\delta,b+\delta]$. Denote by $k$ the integer part of the number $b+\frac{1}{2}$. Then by definition $w(b)=w_k$. If $\delta \leq 1$, then it follows from $\{w_k\}\in \tilde{\Delta}_{2}$ that
\[w(2B)=\int_{b-2\delta}^{b+2\delta}w(x)dx \leq 4\delta (w_{k-1}+w_{k}+w_{k+1}) \leq 4\delta \cdot D \cdot w_{k} \leq 4D \cdot w(B). \]
If $\delta \in [r, r+1)$ for some $r \in \mathbb{N}$, then
\[w(2B)=\int_{b-2\delta}^{b+2\delta}w(x)dx \leq \sum_{m=k-2r}^{k+2r}w_{m} \leq D \cdot \sum_{m=k-r}^{k+r}w_{m} \leq D^{2} \sum_{\{k \in \mathbb{Z}: \ m \in [k-\frac{\delta}{2}, k+\frac{\delta}{2}] \} }w_{k} \leq 2D^{2}w(B). \]
Therefore, $w \in \Delta_2$. Then it follows from Lemma 2.2 that there exists a constant $D_{1}>1$ such that for any interval $B$
\[ w(2B)\geq D_1 w(B).\]
If we take $B=[m-n-\frac{1}{2},m+n+\frac{1}{2}]$, then we have
\[ \sum_{k=m-(2n+1)}^{m+(2n+1)} w_k \ge w(2B) \ge D_1 w(B) = D_1 \sum_{k=m-n}^{m+n} w_k. \]
This completes the proof of Lemma 2.3.

Now, let us prove Lemma 2.4. Let $\{w_k\}\in \tilde{A}_p$, $1<p<\infty$. Let us take the interval $B=(x-r,x+r)$ and denote the integer parts of the numbers $x-r+\frac{1}{2}$ and $x+r+\frac{1}{2}$ by $m$ and $n$. Let
\[ \delta_k=|B \cap[m-\frac{1}{2}+(k-1),m-\frac{1}{2}+k]|, \ \ k=\overline{1,n-m+1}, \ \  \delta=\max_{k=\overline{1,n-m+1}} \delta_k. \]
Then
\[ \left(\frac{1}{|B|}\int_{B} w(x)dx\right) \left(\frac{1}{|B|} \int_{B}w(x)^{-\frac{1}{p-1}}dx \right)^{p-1}\]
\[ = \frac{\left(\sum_{k=1}^{n-m+1} \delta_k w_{m+k-1}\right) \left(\sum_{k=1}^{n-m+1} \delta_k w_{m+k-1}^{-\frac{1}{p-1}} \right)^{p-1}}{\left(\sum_{k=1}^{n-m+1}\delta_k \right)^p} \]
\[ \le \left( \frac{\delta}{\sum_{k=1}^{n-m+1}\delta_k} \right)^p \cdot \left(\sum_{k=m}^{n} w_{k}\right) \left(\sum_{k=m}^{n} w_{k}^{-\frac{1}{p-1}} \right)^{p-1} \]
\[ \le 3^p \cdot \frac{\left(\sum_{k=m}^{n} w_{k}\right) \left(\sum_{k=m}^{n} w_{k}^{-\frac{1}{p-1}} \right)^{p-1}}{(n-m+1)^p} \le 3^p \cdot [w]_{\tilde{A}_p}\]
It follows that $w(x) \in A_p$. Then from Lemma 2.1 we obtain $w(x)\in \Delta_2$. Now we will prove that $ \{w_k\} \in \tilde{\Delta}_2$, that is, inequality (\ref{2.2}) is satisfied. For $\delta < \frac{1}{2}$, inequality (\ref{2.2}) is obvious. For $\delta \ge \frac{1}{2}$, we have
\[ \sum_{ \{k \in \mathbb{Z}: \ k \in [m-2\delta, m+2\delta] \} } w_k \leq w([m-4\delta, m+4\delta]) \]
\[ \leq D^3 w([m-\frac{\delta}{2}, m+\frac{\delta}{2}]) \leq D^3\sum_{ \{k \in \mathbb{Z}: \ k \in [m-\delta, m+\delta] \} } w_k. \]
This completes the proof of Lemma 2.4.
\\  \\
\textbf{Corollary 2.1.} If $\{w_k\} \in \tilde{A}_p$ for some $1 \leq p < \infty$, then there exists $D_1 >1$ such that for any $m \in \mathbb{Z}$ and $i \in \mathbb{N}$
\[ \sum_{k=m-2^i}^{m+2^i} w_k \ge D_1^{i-1} (w_{k-1}+w_k+w_{k+1}). \] 
\section{Boundedness of the Discrete Hilbert Transform on Discrete Weighted Morrey Spaces}

\textbf{Theorem 3.1.} Let $1\le p<\infty$, $0\le \lambda \le \frac{1}{p}$. If $\{w_k\} \in \tilde{A}_p$, then for all $b=\{b_m\}_{m \in \mathbb{Z}} \in m_{\lambda,p,w}$ the sequence $H(b) = \{(Hb)_n\}_{n \in \mathbb{Z}}$ is well defined. 
\\  \\
\textit{Proof.} Let $b=\{b_n\}_{n\in \mathbb{Z}} \in m_{\lambda,p,w}$ , $1<p<\infty$ , $0<\lambda<\frac{1}{p}$. Then for any $n \in \mathbb{N}$
\[ |(Hb)_n| \leq \sum_{m\neq n} \frac{|b_m|}{|n-m|} = \sum_{i=1}^{\infty} \sum_{2^{i-1}\leq|n-m|< 2^i} \frac{|b_m|{w_m}^\frac{1}{p}}{|n-m|{w_m}^\frac{1}{p}} \leq \sum_{i=1}^{\infty} \frac{1}{2^{i-1}} \sum_{2^{i-1}\leq|n-m|< 2^i} \frac{|b_m|{w_m}^\frac{1}{p}}{{w_m}^\frac{1}{p}} \] 
\[ \leq \sum_{i=1}^{\infty} \frac{1}{2^{i-1}} \left(\sum_{|n-m|<2^i} |b_m|^p w_m\right)^\frac{1}{p} \left(\sum_{|n-m|<2^i} {w_m}^{-\frac{1}{p-1}} \right)^\frac{p-1}{p} \] 
\[\leq \sum_{i=1}^{\infty} \frac{1}{2^{i-1}} \left(\sum_{|n-m|<2^i} |b_m|^p w_m\right)^\frac{1}{p}  \frac{[w]_{\tilde{A}_p}^\frac{1}{p} (2^{i+1}-1)}{\left(\sum_{|n-m|<2^i} w_m\right)^{\frac{1}{p}}} \] 
\[ \leq  \sum_{i=1}^{\infty} \frac{1}{2^{i-1}} \frac{\left(\sum_{|n-m|<2^i} |b_m|^p w_m\right)^\frac{1}{p}}{\left(\sum_{|n-m|<2^i} w_m\right)^\lambda} \frac{[w]_{\tilde{A}_p}^\frac{1}{p} (2^{i+1}-1)}{\left(\sum_{|n-m|<2^i} w_m\right)^{\frac{1}{p}-\lambda}} \]
\[ \leq 4 ||b||_{m_{\lambda,p,w}} [w]_{\tilde{A}_p}^\frac{1}{p} \sum_{i=1}^{\infty} \left( \sum_{|n-m|<2^i} w_m\right)^{\lambda-\frac{1}{p}}. \]
Denoting
\[ \delta=\log_2 D_1 >0, \]
and using Corollary 2.1, we have
\[ \sum_{m=n-2^i}^{n+2^i} w_m \geq D_1 ^{i-1} \sum_{m=n-1}^{n+1} w_m =2^{(i-1)\delta} (w_{n-1}+w_n+w_{n+1}). \]
Then
\[ |(Hb)_n|\leq 4||b||_{m_{\lambda,p,w}} [w]_{\tilde{A}_p}^\frac{1}{p} \frac{2^{\delta(\frac{1}{p}-\lambda)}}{\left(2^{\delta(\frac{1}{p}-\lambda)}-1\right)(w_{n-1}+w_n+w_{n+1})^{\frac{1}{p}-\lambda}} <\infty. \]
Therefore, the Hilbert transform of the sequence $b$ is well defined. This completes the proof of the Theorem 3.1.
\\  \\
\textbf{Theorem 3.2.} Let $1\le p<\infty$, $0\le \lambda \le \frac{1}{p}$. If $\{w_k\} \in \tilde{A}_p$, then for any $b=\{b_m\}_{m \in \mathbb{Z}} \in m_{\lambda,p,w}$ we have $H(b) \in m_{\lambda,p,w}$ and there exists $C_{\lambda,p,w} >0$, such that 
\[ \|H(b)\|_{m_{\lambda,p,w}} \leq C_{\lambda,p,w} \|b\|_{m_{\lambda,p,w}}\] 
for any $b\in m_{\lambda,p,w}$.
\\ \\ 
\textit{Proof.} We define $f(x)$ as $b_k$ for $x\in [k-\frac{1}{4},k+\frac{1}{4}]$, $ k\in \mathbb{Z}$, and $0$ elsewhere and $w(x)$ as $w_k$ for $x\in[k-\frac{1}{4},k+\frac{1}{4}]$ and linearly between. It is easy to verify if $\{w_k\}\in \tilde{A}_p$ then $w(x)\in A_p$ (see \cite{Hunt}). We are going to show that $f\in M_{\lambda,p,w}$. Let us take the interval $B=(x-r,x+r)$ and denote the integer parts of the numbers $x-r+\frac{1}{2}$ and $x+r+\frac{1}{2}$ by $m$ and $n$. Let
\[ \delta_k=|B \cap[m-\frac{1}{4}+(k-1),m-\frac{1}{4}+k]|, \ \ k=\overline{1,n-m+1}, \ \  \delta=\max_{k=\overline{1,n-m+1}} \delta_k. \]
Then
\[ \frac{\left(\int_{B}|f(t)|^p w(t)dt\right)^\frac{1}{p}}{\left(\int_{B} w(t)dt\right)^\lambda} = \left(\int_{B} w(t)dt\right)^{-\lambda} \cdot \left(\sum_{k=1}^{n-m+1} \delta_k w_{m+k-1}|b_{m+k-1}|^p \right)^{\frac{1}{p}} \]
\[ \le \left( \int_{B} w(t)dt \right)^{-\lambda} \cdot \delta^{\frac{1}{p}} \cdot \left(\sum_{k=m}^{n} w_{k}|b_{k}|^p \right)^{\frac{1}{p}}. \]
If $\delta=0$, then
\[ \frac{\left(\int_{B}|f(t)|^p w(t)dt\right)^\frac{1}{p}}{\left(\int_{B} w(t)dt\right)^\lambda} = 0.\]
Let us $\delta>0$. If $n=m$, that is, if $B \subset [m-\frac{1}{2}, m+\frac{1}{2})$ for some $m \in \mathbb{Z} $, then
\[ \left(\int_{B} w(t)dt\right) \ge \delta_1 \cdot w_m = \delta \cdot w_m=\delta \cdot \sum_{k=m}^n w_k ;\]
if $n=m+1$, that is, if $[m-\frac{1}{2}, m+\frac{1}{2}] \subset B \subset [m-\frac{1}{2}, m+\frac{3}{2})$ for some $m \in \mathbb{Z} $, then it follows from Lemma 2.4 that
\[\left(\int_{B} w(t)dt\right) \ge \delta_1 \cdot w_m+\delta_2 \cdot w_{m+1} \ge \delta \cdot \min \{w_m, w_{m+1} \}\ge \frac{\delta}{D}\cdot \sum_{k=m}^n w_k ;\]
if $n \ge m+2$, that is, if $[m-\frac{1}{2}, m+\frac{3}{2}] \subset B \subset [m-\frac{1}{2}, n+\frac{1}{2})$, then it follows from Lemma 2.4 that
\[\left(\int_{B} w(t)dt\right) \ge \delta \sum_{k=m+1}^{n-1} w_k \ge \frac{\delta}{D} \cdot \sum_{k=m}^n w_k .\]
Therefore,
\[ \frac{\left(\int_{B}|f(t)|^p w(t)dt\right)^\frac{1}{p}}{\left(\int_{B} w(t)dt\right)^\lambda} \leq \delta^{\frac{1}{p}-\lambda} \cdot D^{\lambda}\cdot\frac{\left(\sum_{k=m}^{n} w_{k}|b_{k}|^p \right)^{\frac{1}{p}}}{\left(\sum_{k=m}^{n} w_{k}\right)^{\lambda}} \]
It follows that
\[ \|f\|_{M_{\lambda,p,w}} \leq \delta^{\frac{1}{p}-\lambda} \cdot D^{\lambda} \|b\|_{m_{\lambda,p,w}}. \]

We denote \[ (Tb)_n =\sup_{k>0}\left|\sum_{m: \ |n-m|\geq k}\frac{b_m}{n-m}\right|, \indent \indent (Mf)(x)=\sup_{y\neq x}\frac{1}{y-x} \int_{x}^{y}|f(t)|dt \] 
\[ S(f)(x)=\sup_{\epsilon>0}\left|\int_{|x-y|>\epsilon}\frac{f(y)}{x-y}dy \right|. \]
In \cite{Hunt}, it was shown that for any $x\in [j-\frac{1}{4},j+\frac{1}{4}]$
\[ \left| \int_{|x-y|>n+\frac{1}{2}} \frac{2f(y)}{x-y}dy - \sum_{k: \ |k-j|>n} \frac{b_k}{j-k} \right| < \int_{\frac{1}{2}<|x-y|} \frac{|f(y)|}{(x-y)^2}dy. \]
It follows that
\begin{equation}\label{2.4}
(Tb)_j \leq 2S(f)(x) + 2\int_{\frac{1}{2}\leq|x-y|}\frac{|f(y)|}{(x-y)^2}dy 
\end{equation}
for $x\in [j-\frac{1}{4},j+\frac{1}{4}]$.
By integrating by parts, it is clear that the second term on the right side of (4) is bounded by $4(Mf)(x)$. Therefore for any $x\in [j-\frac{1}{4},j+\frac{1}{4}]$
\begin{equation}\label{2.5}
(Tb)_j \leq 2S(f)(x) + 4 (Mf)(x).
\end{equation}
In \cite{Komori}, it was shown that if $w(x)\in A_p$, then the Hardy-Littlewood maximal operator $M$ is bounded in weighted Morrey space $M_{\lambda,p,w}$. Similar to the proof of (\cite{Komori}, Theorem 3.3), we prove the boundedness of the operator $S$ in the space $M_{\lambda,p,w}$. Fix an interval $B=(x_0-r,x_0+r)$, and decompose $f=f_1+f_2$ with $f_1=f\cdot \chi_{2B}$. Then
\[ \int_{B}|Sf(x)|^pw(x)dx \leq \int_{B}|Sf_1(x)|^pw(x)dx+\int_{B}|Sf_2(x)|^pw(x)dx. \]
Using the fact that if $w(x)\in A_p$, then $S$ is bounded on $L_p(w)$ (see \cite{Hunt}), we can get
\[ \int_{B}|Sf_1(x)|^pw(x)dx \leq \int_{R}|Sf_1(x)|^pw(x)dx\]
\[\leq C \int_{R}|f_1(x)|^pw(x)dx=\int_{B}|f(x)|^pw(x)dx. \]
Hence we have
\begin{equation}\label{2.6}
\int_B |Sf_1(x)|^pw(x)dx \leq C \cdot \|f\|_{M_{\lambda,p,w}}^p \cdot w(B)^{\lambda p}
\end{equation}
For $x \in B$ and $y \notin 2B$ we have $|x-y| \ge \frac{1}{2} |x_0-y|$. Thus, for $x \in B$ we get
\begin{equation}\label{2.7}
|Sf_2(x)| \leq \int_{R}\frac{|f_2(y)|}{|x-y|}dy \leq C \cdot \int_{|x_0-y| \ge 2r}\frac{|f(y)|}{|x_0-y|}dy.
\end{equation}
It follows from Holder’s inequality and $w \in A_p$ that
\[ \int_{|x_0-y| \ge2r}\frac{|f(y)|}{|x_0-y|}dy =\sum_{j=1}^{\infty} \int_{2^jr \le|x_0-y|<2^{j+1}r}\frac{|f(y)|}{|x_0-y|}dy \]
\[ \leq \sum_{j=1}^{\infty} \frac{1}{2^{j}r} \int_{|x_0-y|<2^{j+1}r}|f(y)|dy \]
\[ \leq \sum_{j=1}^{\infty} \frac{1}{2^{j}r} \left(\int_{|x_0-y|<2^{j+1}r}|f(y)|^pw(y)dy \right)^{\frac{1}{p}} \left(\int_{|x_0-y|<2^{j+1}r}w(y)^{-\frac{1}{p-1}}dy \right)^{\frac{p-1}{p}}  \]
\[ \leq C \cdot \|f\|_{M_{\lambda,p,w}} \cdot \sum_{j=1}^{\infty} \frac{1}{w(2^{j+1}B)^{\frac{1}{p}-\lambda}}  \]
Using (\ref{2.7}), Lemma 2.1 and Lemma 2.2, we get
\[ \int_B |Sf_2(x)|^pw(x)dx \leq C \cdot \|f\|_{M_{\lambda,p,w}}^p \cdot w(B)^{\lambda p} \cdot \sum_{j=1}^{\infty} \left(\frac{w(B)}{w(2^{j+1}B)} \right)^{1-p\lambda} \]
\begin{equation}\label{2.8}
\leq C \cdot \|f\|_{M_{\lambda,p,w}}^p \cdot w(B)^{\lambda p}
\end{equation}
It follows from (\ref{2.6}) and (\ref{2.8}) that the operator $S$ is bounded in the space $M_{\lambda,p,w}$. From (\ref{2.5}) we have
\[ \sum_{k\in S_{m,n}} |Tb_k|^p w_k \leq  \sum_{k\in S_{m,n}} 2\int_{k-\frac{1}{4}}^{k+\frac{1}{4}} |2S(f)+4M(f)|^p w(x)dx \leq 4\int_{m-n-\frac{1}{4}}^{m+n+\frac{1}{4}} |S(f)+2M(f)|^p w(x)dx \]
Since
\[\int_{m-n-\frac{1}{4}}^{m+n+\frac{1}{4}}w(x)dx = \sum_{k=m-n}^{m+n} \int_{k-\frac{1}{4}}^{k+\frac{1}{4}}w(x)dx +\sum_{k=m-n}^{m+n-1}\int_{k+\frac{1}{4}}^{k+\frac{3}{4}}w(x)dx \leq \sum_{k=m-n}^{m+n} w_k, \] 
then
\[ \frac{\left(\sum_{k\in S_{m,n}} |Tb_k|^p w_k\right)^\frac{1}{p}} {\left(\sum_{k\in S_{m,n}} w_k\right)^\lambda} \leq\frac{\left( 4\int_{m-n-\frac{1}{4}}^{m+n+\frac{1}{4}} |S(f)+2M(f)|^p w(x)dx\right)^\frac{1}{p}}{\left(\int_{m-n-\frac{1}{4}}^{m+n+\frac{1}{4}}w(x)dx\right)^\lambda} \].
\[ \leq 4^{\frac{1}{p}} ||S(f)+2M(f)||_{M_{\lambda,p,w}}\leq 4^{\frac{1}{p}} \left(||S(f)||_{M_{\lambda,p,w}} + 2||M(f)||_{M_{\lambda,p,w}} \right). \]
Therefore, $T(b) \in m_{\lambda,p,w}$ and
\[\|T(b)\|_{m_{\lambda,p,w}} \leq  4^{\frac{1}{p}}\left(||S(f)||_{M_{\lambda,p,w}} + 2||M(f)||_{M_{\lambda,p,w}} \right)  \|b\|_{m_{\lambda,p,w}} \]
\[\leq  4^{\frac{1}{p}}\left(||S||_{M_{\lambda,p,w} \to M_{\lambda,p,w}} + 2||M||_{M_{\lambda,p,w} \to M_{\lambda,p,w}} \right)  \|f\|_{M_{\lambda,p,w}} \]
\[\leq  4^{\frac{1}{p}}\left(||S||_{M_{\lambda,p,w} \to M_{\lambda,p,w}} + 2||M||_{M_{\lambda,p,w} \to M_{\lambda,p,w}} \right) \delta^{\frac{1}{p}-\lambda} \cdot D^{\lambda}  \|b\|_{m_{\lambda,p,w}} \]
 On the other side, it follows from inequality
\[ |(Hb)_n|\leq |(Tb)_n|\]
that $ H(b) \in m_{\lambda,p,w}$ and
\[ \|H(b)\|_{m_{\lambda,p,w}} \leq 4^{\frac{1}{p}}\left(||S||_{M_{\lambda,p,w} \to M_{\lambda,p,w}} + 2||M||_{M_{\lambda,p,w} \to M_{\lambda,p,w}} \right) \delta^{\frac{1}{p}-\lambda} \cdot D^{\lambda}  \|b\|_{m_{\lambda,p,w}}\] 
for any $b\in m_{\lambda,p,w}$. This completes the proof of the Theorem 3.2.

\newpage
\renewcommand{\refname}

\end{document}